\documentclass[a4paper,12pt,reqno]{amsart}
\usepackage{latexsym}
\usepackage{amsmath,amsthm,amssymb,amscd}
\usepackage{fullpage}
\usepackage{xcolor}
\usepackage{stmaryrd}
\usepackage{parskip}
\usepackage{mathtools}
\mathtoolsset{showonlyrefs}

\usepackage[activate={true,nocompatibility},final,tracking=true,kerning=true,spacing=true]{microtype}
\usepackage[pdfencoding=unicode,pdftex,bookmarks=true,bookmarksnumbered]{hyperref}
\definecolor{citecolour}{rgb}{0.0, 0.0, 0.8}
\colorlet{linkcolour}{green!50!black}
\hypersetup{colorlinks,breaklinks,
		linkcolor=linkcolour,citecolor=citecolour,
		filecolor=linkcolour, urlcolor=linkcolour}
\usepackage{version}

\newtheorem{prevtheorem}{Theorem}

\newtheorem{theorem}{Theorem}[section]
\newtheorem{lemma}[theorem]{Lemma}

\newtheorem{remark}[theorem]{Remark}

\numberwithin{equation}{section}
\usepackage[shortlabels]{enumitem}
\begingroup
    \makeatletter
    \@for\theoremstyle:=definition,remark,plain\do{%
        \expandafter\g@addto@macro\csname th@\theoremstyle\endcsname{%
            \addtolength\thm@preskip\parskip
            }%
        }

\newcommand{\R}{\mathbb{\R}}

\newcommand{\s}{\mathfrak{s}}

\includeversion{comment}
\renewcommand{\leq}{\leqslant}
\renewcommand{\geq}{\geqslant}

\newenvironment{proofof}{{\bf {Proof.} }}{\hfill $\blacksquare$ \\}

\begin{document}
\title{Finite non-cyclic $p$-groups whose number of subgroups is minimal}
\author[S. Aivazidis and 
T.\,W. M\"uller]{Stefanos Aivazidis$^{\dagger}$ and
Thomas M\"uller$^*$}
\address{$^{\dagger*}$Stockholm, Sweden.}
\email{stefanosaivazidis@gmail.com}
\address{$^*$School of Mathematical Sciences, Queen Mary
\& Westfield College, University of London,
Mile End Road, London E1 4NS, United Kingdom.}
\email{t.w.muller@qmul.ac.uk}
\begin{abstract}
Recent results of Qu and T\u{a}rn\u{a}uceanu describe explicitly the finite $p$-groups which are not elementary abelian and have the property that the number of their subgroups is maximal among $p$-groups of a given order. We complement these results from the bottom level up by determining completely the non-cyclic finite $p$-groups whose number of subgroups among $p$-groups of a given order is minimal. 
\end{abstract}
\maketitle
\section{Introduction}\label{Sec:Intro}
Let 
\begin{equation}
M_p(1,1,1) \coloneqq \left\langle a,b, c \mid a^p = b^p = c^p = 1, [a,b] = c, [c,a] = [c, b] = 1 \right\rangle. 
\end{equation}
Write $\mathfrak{s}_k(G)$ for the set of subgroups of order $p^k$ of the $p$-group $G$ and 
$\mathfrak{s}(G)$ for the set of all subgroups of $G$. Recently, Qu obtained the following result.

\begin{theorem}[{\cite[Thm.~1.4]{Qu}}]
\label{Thm:Qu}
Let $G$ be a group of order $p^n$, where $p$ is an odd prime, 
and $\widetilde{G}  = M_p(1, 1, 1) \times C_p^{n-3}$. If G is not elementary abelian, then for all $k$ such that $1 \leq k \leq n$, we have 
$\left\lvert \mathfrak{s}_k(G)\right\rvert \leq \lvert \mathfrak{s}_k(\widetilde{G})\rvert$. 
In particular, if $2 \leq k \leq n-2$, then $\lvert \mathfrak{s}_k(G)\rvert < \lvert\s_k(\widetilde{G})\rvert$.
\end{theorem}

Thus Qu finds the $p$-groups whose number of subgroups of possible order is maximal 
except for elementary abelian $p$-groups when $p>2$. The analogue of Qu's result for $p=2$
was recently obtained by T\u{a}rn\u{a}uceanu (cf.~\cite{tar}). We wish to complement Qu's and 
T\u{a}rn\u{a}uceanu's results from the bottom level up. In particular, 
in Theorems~\ref{Thm:p>2} and~\ref{Thm:p=2} we determine completely the structure of those
finite non-cyclic $p$-groups whose number of subgroups is minimal.
\section{Finite $p$-groups realising the second minimal level of subgroup numbers}

We begin with the following preparatory lemma.

\begin{lemma}\label{Lem:Butler}
Suppose that $G$ is an abelian $p$-group of order $p^{\lambda}$, and that $G$ has $p+1$ subgroups 
of order $p^r$ for all $r$ such that $1 \leq r \leq \lambda-1$. Then $G \cong C_{p^{\lambda-1}} \times C_p$.
\end{lemma}

\begin{proofof}
By a well-known result, the number of subgroups of type $\nu = (\nu_1,\ldots,\nu_\ell)$ 
in an abelian $p$-group of type $\mu = (\mu_1, \ldots, \mu_m)$ is given by the formula
\begin{equation}\label{Eq:psubCount}
\prod_{i\geq 1} p^{\nu_{i+1}'(\mu_i' - \nu_i')} {\mu_i' - \nu_{i+1}' \brack \nu_i' - \nu_{i+1}'}_p,
\end{equation}
where $\mu', \nu'$ are the conjugates of the partitions $\mu$ and $\nu$, respectively, and
\begin{equation}
{n \brack k}_p = \prod_{i=0}^{k-1} \frac{1 - p^{n-i}}{1-p^{i+1}}
\end{equation}
is the number of $k$-dimensional subspaces of an $n$-dimensional vector space over the field $\mathbb{Z}/p\mathbb{Z}$; 
see, for instance,~\cite[Eqn.~(1)]{Butler}. 
Since $G$ is non-cyclic of type $\mu = (\mu_1, \ldots, \mu_m)$ say, it follows that $m\geq 2$ (and $\mu_2 >0$), thus $\mu_1' \geq 2$. Now, let us the count the number of subgroups of type $\nu = (1,0,\ldots,0)$. Observe that $\nu'=(1,0,\ldots,0)$, so according to formula \eqref{Eq:psubCount} there are 
\[
r \coloneqq {\mu_1' \brack 1}_p \cdot \prod_{i\geq 2} p^{\nu_{i+1}'(\mu_i' - \nu_i')} {\mu_i' - \nu_{i+1}' \brack \nu_i' - \nu_{i+1}'}_p,
\]
such subgroups. However, $r \leq p+1$ since by assumption the number of subgroups of order $p$ in $G$ is $p+1$, and $r \geq p+1$, owing to
\[
{\mu_1' \brack 1}_p \geq {2 \brack 1}_p = p+1.
\]
Thus $\mu_1'=2$, which implies that $m=2$ and that $G$ is of type $\mu=(\mu_1,\mu_2)$. It follows that 
\[
\mu'=(\underbrace{2,2,\ldots,2}_{\mu_2},\underbrace{1,1,\ldots,1}_{\mu_1 -\mu_2})
\]
Suppose that $\mu_2 \geq 2$. Next, we count subgroups of type $\nu=(\mu_1-1,1)$ in $G$. Note that 
\[
\nu'=(2,\underbrace{1,\ldots,1}_{\mu_1-1}).
\]
Moreover, the second term in the product formula for $\nu=(\mu_1-1,1)$ evaluates to $p$, thus the number of subgroups of type $\nu=(\mu_1-1,1)$ in $G$ is at least $p$. Further, the $(\mu_2-1)$-term in the product formula for $\nu=(\mu_1,0)$ evaluates to $p$ as well, and thus we deduce that the number of subgroups of type $\nu=(\mu_1,0)$ in $G$ is, again, at least $p$.
Therefore, we get at least $2p>p+1$ subgroups of order $p^{\mu_1}$ in $G$, contrary to our assumption that the number of subgroups of order $p^{\mu_1}$ is $p+1$.
This contradiction shows that $\mu_2=1$ and $\mu_1=\lambda-1$, so $G = C_{p^{\lambda-1}} \times C_p$ is of the type we asserted.
\end{proofof}

We shall present our main result as two separate theorems, 
dealing with the cases $p>2$ and $p=2$ respectively,
since the $p=2$ case, although ultimately similar to the $p>2$ case, 
presents a somewhat erratic behaviour at small values. 

\begin{prevtheorem}\label{Thm:p>2}
Let $p$ be an odd prime, $G$ a $p$-group of order $p^{\lambda}$. If $G$ is not the cyclic group 
$C_{p^{\lambda}}$, then $\lvert \mathfrak{s}(G)\rvert \geq (p+1)(\lambda -1) +2$, with equality if and only if 
$G \cong C_{p^{\lambda-1}} \times C_p$ or 
\begin{equation}\label{Eq:modpgroup}
G \cong M_{p^\lambda} \coloneqq \left\langle a,b \mid a^{p^{\lambda-1}}=b^p=1, a^b=a^{1+p^{\lambda-2}}\right\rangle.
\end{equation}
\end{prevtheorem}

\begin{proofof}
Given a non-cyclic $p$-group $G$, $p$ odd, a well-known theorem due to Kulakoff~\cite[Satz 1]{Kulakoff} asserts that 
\[\lvert \s_k(G)\rvert \equiv p+1\;(\bmod\; p^2)\]
for all $k$ such that $1 \leq k \leq \lambda-1$. Thus, in particular, $\lvert \s_k(G)\rvert \geq p+1$, and therefore
\[
\lvert \mathfrak{s}(G)\rvert = 
\sum_{k=0}^{\lambda}\lvert \mathfrak{s}_k(G)\rvert = 
2 + \sum_{k=1}^{\lambda-1}\lvert\s_k(G)\rvert \geq 2 + (\lambda-1)(p+1).
\]
This proves the first part of our assertion.

Now, we assume (as we may) that $\lambda \geq 3$. Applying a result of 
Lindenberg \cite[Folgerung~3.4]{Lindenberg}, we have that if $G = C_{p^{\lambda-1}} \rtimes C_p$ (the implied 
action of $C_p$ on $C_{p^{\lambda-1}}$ may well be trivial; we only require that $G$ be a split extension) then 
$G$ has, apart from the trivial subgroup and the whole group $G$, $p+1$ subgroups of order $p^j$, 
$1 \leq j \leq \lambda-1$, and thus 
$\lvert \mathfrak{s}(C_{p^{\lambda-1}} \rtimes C_p)\rvert = (p+1)(\lambda-1)+2$ 
for any such split extension.

Our goal now is to establish that if a (necessarily non-cyclic) $p$-group $G$ of order $|G|=p^\lambda$ has $p+1$ subgroups 
of order $p^j$, for all $j$ such that $1 \leq j \leq \lambda-1$, then $G$ is a split extension $C_{p^{\lambda-1}} \rtimes C_p$.

If $G$ is abelian, then the claim follows from Lemma~\ref{Lem:Butler}. Now, assume that $G$ is non-abelian. Since $G$ is a 
$p$-group, $p>2$, this implies that $G$ is not a Dedekind group. Of the $p+1$ subgroups of order $p^j$, for each $j$ such 
that $1\leq j \leq \lambda-1$, one is certainly normal, since $p$-groups have normal subgroups of each possible order. The 
other $p$ are either all normal, or lie in the same conjugacy class. It follows that the non-normal subgroups of $G$ for each 
possible order are all conjugate, and thus $G$ is a CO-group; cf. top of~\cite[Sec. 58]{bj2}. Janko's theorem 
(see~\cite[Thm.~58.3]{bj2}) now yields 
$G \cong M_{p^\lambda}$, where $M_{p^\lambda}$ was defined in \eqref{Eq:modpgroup}.
Our proof is complete.
\end{proofof}

\begin{remark}
We note here that there may, in principle, exist many non-isomorphic semidirect products $C_{p^{\lambda-1}} \rtimes C_p$ which 
are not direct products. However, an early theorem due to Burnside asserts that the only non-abelian $p$-group, $p$ odd, which 
has a cyclic maximal subgroup is 
$M_{p^\lambda}$; see \cite[Chap.~VIII, Sec. 109]{Burnside}.
\end{remark}

Recall that the generalised quaternion group of order $2^{\lambda}$ is the group defined by the presentation
\[
Q_{2^{\lambda}} \coloneqq \left\langle x,y \mid x^{2^{\lambda-1}} = 1, y^2 = x^{2^{\lambda-2}}, x^y = x^{-1}\right\rangle.
\]

Next, we address the $p=2$ case. 

\begin{prevtheorem}\label{Thm:p=2}
Let $G$ be a $2$-group of order $2^{\lambda}$. If $G$ is not the cyclic group 
$C_{2^{\lambda}}$ and $\lambda \geq 5$, then $\lvert\mathfrak{s}(G)\rvert \geq 3\lambda -1$, with equality if and only if 
$G \cong C_{2^{\lambda-1}} \times C_2$ or 
$G \cong M_{2^{\lambda}}$, where $M_{2^{\lambda}}$ is as in \eqref{Eq:modpgroup} with $p=2$.
If $\lambda = 3$ and $G$ is not $C_8$, then $\lvert\mathfrak{s}(G)\rvert \geq 6$ with equality if and only if $G \cong Q_8$, while if 
$\lambda = 4$ and $G$ is not $C_{16}$, then 
$\lvert\mathfrak{s}(G)\rvert \geq  11$, with equality if and only if $G \cong Q_{16}$, or $G \cong M_{16}$, or $G \cong C_8 \times C_2$.
\end{prevtheorem}

\begin{proofof}
There are 5 groups of order $8=2^3$, and 14 groups of order $16=2^4$. We use GAP~\cite{GAP4} to obtain a full list of the isomorphism 
classes of groups in each case, and ask GAP for the total number of subgroups of each group in the list. Our claim for $\lambda=3$ and 
$\lambda=4$ is now a simple matter of inspection.

We may thus assume that $\lambda \geq 5$. By Frobenius' generalisation of Sylow's theorem, we have
\[
\lvert\s_{\mu}(G)\rvert \equiv 1\;(\bmod\; 2),
\]
for all $\mu$ such that $0 \leq \mu \leq \lambda$. 

Moreover, since $G \ncong C_{2^{\lambda}}$, we have 
$\lvert\s_{\mu}(G)\rvert \neq 1$, $1 \leq \mu \leq \lambda-1$
by \cite[Prop.~1.3]{bj1} unless $G \cong Q_{2^{\lambda}}$ and $\mu=1$.
It is well-known that $Q_{2^{\lambda}}$ has a unique involution $t$, which generates its centre, and affords the quotient 
$Q_{2^{\lambda}}\big/\langle t \rangle \cong D_{2^{\lambda -1}}$.
Since $\langle t \rangle$ is the unique subgroup of order 2 in $Q_{2^{\lambda}}$, it follows easily that 
\[
\lvert\mathfrak{s}(Q_{2^{\lambda}})\rvert = \lvert\mathfrak{s}(D_{2^{\lambda-1}})\rvert + 1 = \tau(2^{\lambda-2}) + \sigma(2^{\lambda-2}) + 1 = 2^{\lambda-1} + \lambda -1,
\]
where $\lvert\mathfrak{s}(D_{2m})\rvert = \tau(m) + \sigma(m)$ follows from~\cite[Ex.~1]{calhoun}.
Moreover, notice that 
$\lvert\mathfrak{s}(Q_{2^{\lambda}})\rvert \geq 3\lambda -1$ for all $\lambda \geq 4$, with equality precisely when $\lambda=4$.
Since $\lambda \geq 5$, it follows that $G \ncong Q_{2^{\lambda}}$. 
This shows that $\lvert\mathfrak{s}(G)\rvert \geq 3\lambda-1$ for all $\lambda \geq 5$,
and by Lemma~\ref{Lem:Butler} the only abelian group of order $2^{\lambda}$ realising this bound is 
$C_{2^{\lambda-1}} \times C_2$.

It is well known that a 2-group (of order $2^{\lambda}$, $\lambda \geq 5$) is Dedekind non-abelian if and only if 
$G \cong Q_8 \times C_2^{\lambda-3}$. Note here that every group of this type has 
$G_0 \coloneqq Q_8 \times (C_2 \times C_2)$ as a direct factor,
and that $\lvert\s_{2}(Q_8 \times C_2^2)\rvert = 7$. 
Hence, a non-abelian group $G$ realising the bound $\lvert\mathfrak{s}(G)\rvert = 3\lambda-1$ cannot be Dedekind.
By Janko's theorem mentioned previously, the only possible non-abelian group realising the bound 
$\lvert\mathfrak{s}(G)\rvert = 3\lambda-1$ is $M_{2^{\lambda}}$.
But Lindenberg's result applies in the case of $M_{2^{\lambda}}$, since $1+2^{\lambda-2} \not\equiv -1 \;(\bmod\; 4)$. We conclude that the
only groups of order $2^{\lambda}$, $\lambda \geq 5$, which realise the bound 
$\lvert\mathfrak{s}(G)\rvert = 3\lambda-1$ are $C_{2^{\lambda-1}} \times C_2$ and 
$M_{2^{\lambda}}$. This completes our proof.
\end{proofof}

%\bibliographystyle{abbrv}
%\bibliography{Bibliography}

\begin{thebibliography}{10}

\bibitem{bj1}
Y.~Berkovich.
\newblock {\em Groups of prime-power order. {V}ol. 1}, volume~46 of {\em de
  Gruyter Expositions in Mathematics}.
\newblock Walter de Gruyter GmbH \& Co. KG, Berlin, 2008.
\newblock With a foreword by Zvonimir Janko.

\bibitem{bj2}
Y.~Berkovich and Z.~Janko.
\newblock {\em Groups of prime-power order. {V}ol. 2}, volume~47 of {\em de
  Gruyter Expositions in Mathematics}.
\newblock Walter de Gruyter GmbH \& Co. KG, Berlin, 2008.

\bibitem{Burnside}
W.~Burnside.
\newblock On some properties of groups whose orders are powers of primes.
\newblock {\em Proc. London Math. Soc. (2)}, 11:225--245, 1913.

\bibitem{Butler}
L.~M. Butler.
\newblock A unimodality result in the enumeration of subgroups of a finite
  abelian group.
\newblock {\em Proc. Amer. Math. Soc.}, 101(4):771--775, 1987.

\bibitem{calhoun}
W.~C. Calhoun.
\newblock Counting the subgroups of some finite groups.
\newblock {\em Amer. Math. Monthly}, 94(1):54--59, 1987.

\bibitem{GAP4}
The GAP~Group,
\newblock \emph{GAP -- Groups, Algorithms, and Programming, Version 4.10.1}; 
2019,   \verb+(https://www.gap-system.org)+.

\bibitem{Kulakoff}
A.~Kulakoff.
\newblock \"{U}ber die {A}nzahl der eigentlichen {U}ntergruppen und der
  {E}lemente von gegebener {O}rdnung in {$p$}-{G}ruppen.
\newblock {\em Math. Ann.}, 104(1):778--793, 1931.

\bibitem{Lindenberg}
W.~Lindenberg.
\newblock \"{U}ber die {S}truktur zerfallender bizyklischer {$p$}-{G}ruppen.
\newblock {\em J. Reine Angew. Math.}, 241:118--146, 1970.

\bibitem{Qu}
H.~Qu.
\newblock Finite non-elementary abelian {$p$}-groups whose number of subgroups
  is maximal.
\newblock {\em Israel J. Math.}, 195(2):773--781, 2013.

\bibitem{tar}
M.~{T\u{a}rn\u{a}uceanu}.
\newblock On a conjecture by {H}aipeng {Q}u.
\newblock {\em arXiv e-prints}, page arXiv:1811.07478, Nov 2018.

\end{thebibliography}

\end{document}